\documentclass{amsart}
\usepackage{amsmath,amssymb}
\newtheorem{theorem}{Theorem}
\newtheorem{lemma}[theorem]{Lemma}
\newtheorem{proposition}[theorem]{Proposition}
\newtheorem{corollary}[theorem]{Corollary}

\theoremstyle{definition}

\theoremstyle{remark}
\newtheorem{remark}[theorem]{Remark}

\newtheorem{conjecture}[theorem]{Conjecture}

\def\H{{\mathcal H}}
\def\Z{{\mathbb Z}}
\def\Rad{{\rm Rad}}

\begin{document}
%\frontmatter
\title[Uno's conjecture]
{Uno's conjecture on representation types of Hecke algebras}

\author{Susumu Ariki}
\address{RIMS, Kyoto university,
Kyoto 606--8502, Japan}
\curraddr{}
\email{ariki@kurims.kyoto-u.ac.jp}
\thanks{}

% \date is required for Memoirs; it is the date received by the editor.
\date{}


% Keywords will be recognized for Memoirs only.

\begin{abstract}
Based on a recent result of the author and A.Mathas, we prove that 
Uno's conjecture on representation types of Hecke algebras is true 
for all Hecke algebras of classical type. 
\end{abstract}

\maketitle

\setcounter{page}{1}
%\include{preface}

%\mainmatter
%\include{chap1}

\section{Introduction}

Let $K$ be a field of characteristic $l$, 
$A$ a finite dimensional $K$--algebra. 
We always assume that $K$ is a splitting field of $A$. 
We say that $A$ is of {\sf finite representation type} if 
there are only finitely many 
isomorphism classes of indecomposable $A$--modules. 
If $A$ is a group algebra, then the following theorem 
answers when $A$ is of finite representation type. 

\begin{theorem}{\cite{Hig}\cite{BD}}
Let $G$ be a finite group, $A=KG$ its group algebra. Then 
$A$ is of finite representation type if and only if the 
Sylow $l$--subgroups of $G$ are cyclic. 
\end{theorem}

We restrict ourselves to the case where $G$ is a finite Weyl group 
and see consequences of this result. 

\begin{theorem}
\label{q=1}
Let $W$ be a finite Weyl group. Then $KW$ is of finite 
representation type if and only if $l^2$ does not divide 
the order $|W|$. 
\end{theorem}

For the proof see Appendix. 
Note that Theorem \ref{q=1} does not hold if we consider 
finite Coxeter groups. Dihedral groups $W(I_2(m))$ with 
odd $l$ and $l^2\,|\,m$ are obvious counterexamples. 

Uno conjectured a $q$--analogue of Theorem \ref{q=1}. 
Let $q\in K$ be an invertible 
element, $(W,S)$ a finite Weyl group. We denote by 
$\H_q(W)$ the (one parameter) Hecke algebra associated 
to $W$. The quadratic relation we choose is $(T_s+1)(T_s-q)=0$ 
$(s\in S)$. The {\sf Poincare polynomial} $P_W(x)$ is defined by 
\begin{equation*}
P_W(x)=\sum_{w\in W}x^{l(w)}. 
\end{equation*}

Let $e$ be the smallest positive integer such that 
$q^{e-1}+\cdots+1=0$. If $q=1$ then $e=l$, and if $q\ne 1$ then 
$e$ is the multiplicative order of $q$. As 
the condition that $l^2$ does not divide $|W|$ is the same as 
$(\frac{x^e-1}{x-1})^2$ evaluated at $x=1$ does not divide $P_W(1)$, 
the following is a reasonable guess. We call this 
Uno's conjecture. 

\begin{conjecture} 
Assume that $q\ne 1$ and $K$ a splitting field of $\H_q(W)$. Then 
$\H_q(W)$ is of finite representation type and not semisimple 
if and only if $q$ is a simple root of $P_W(x)=0$, that is, if and 
only if $(x-q)^2$ does not divide $P_W(x)$. 
\end{conjecture}

It is well--known that $\H_q(W)$ is semisimple if and only if 
$P_W(q)\ne 0$. 
See \cite[Proposition 2.3]{G} for example. 
This fact will be used later. 

\begin{remark}
In \cite{U}, it is proved that the conjecture is true for 
$\H_q(I_2(m))$. So, unlike the case of $q=1$, we may ask 
the same question for finite 
Coxeter groups instead of finite Weyl groups. 
\end{remark}

The following theorem is proved in \cite{U}. 

\begin{theorem}{\cite[Proposition 3.7,Theorem 3.8]{U}}
\label{Uno}
Suppose that $q\ne 1$ and denote its multiplicative order by $e$. 
Then $\H_q(A_{n-1})$ is of
finite representation type if and only if $n<2e$.
\end{theorem}

As $P_W(x)=\prod_{i=1}^n \frac{x^i-1}{x-1}$ in this case, 
a primitive $e^{th}$ root of unity is a simple root if and 
only if $n<2e$. 
In particular, the conjecture is true if $W=W(A_{n-1})$. 
The purpose of this article is to prove. 

\begin{theorem}{\bf (Main Theorem)}
\label{main theorem}
Assume that $W$ is of classical type and that 
$K$ is a splitting field of $\H_q(W)$. Then 
$\H_q(W)$ is of finite representation type and not semisimple 
if and only if $q$ is a simple root of $P_W(x)=0$. 
\end{theorem}

We remark that the exceptional cases are settled recently by 
Miyachi \cite{M} under the assumption that 
the characteristic $l$ of the base field $K$ is not too small. 

\section{Reduction to Hecke algebras associated to irreducible Weyl groups}

This is proved by using the complexity of modules. 
Let $A$ be a self--injective finite dimensional $K$--algebra, 
$M$ a finite dimensional $A$--module, $P^{\centerdot}\rightarrow M$ be its 
minimal projective resolution. Then 
the {\sf complexity} $c_A(M)$ is the smallest integer $s\ge 0$ 
such that ${\rm dim}_K(P^t)/(t+1)^{s-1}$ $(t=0,1,\dots)$ is bounded. 
The following lemma is fundamental. 

\begin{lemma}
\label{complexity}
Suppose that $A$ is self--injective as above. Then 

\begin{itemize}
\item[(1)] 
$c_A(M)=0$ if and only if $M$ is a projective $A$--module. 
\item[(2)] 
$A$ is semisimple if and only if $c_A(M)=0$ for all indecomposable 
$A$--modules $M$. 
\item[(3)]
If $A$ is of finite representation type and not semisimple 
then $c_A(M)\le 1$ for all indecomposable $A$--modules $M$ and 
the equality holds for some $M$. 
\end{itemize}
\end{lemma}

\begin{proposition}
\label{reduction}
Let ${\mathcal S}$ be a set of irreducible finite Weyl groups. 
If Uno's conjecture is true for all $\H_q(W)$ with $W\in{\mathcal S}$ 
then the conjecture is true for $\H_q(W_1\times\cdots\times W_r)$ 
with $W_1,\dots,W_r\in{\mathcal S}$. 
\end{proposition}
\begin{proof}
Write $W=W_1\times\cdots\times W_r$. Then 
$\H_q(W)=\H_q(W_1)\otimes\cdots\otimes\H_q(W_r)$ and 
$P_W(x)=P_{W_1}(x)\cdots P_{W_r}(x)$. 

First assume that $q\ne 1$ is a simple root of $P_W(x)=0$. Then 
$q$ is a simple root of $P_{W_i}(x)=0$ and $P_{W_j}(q)\ne 0$ for 
all $j\ne i$. Then $\H_q(W_j)$ for $j\ne i$ are all semisimple 
and $\H_q(W_i)$ is of finite representation type and not semisimple. 
Thus $\H_q(W)$ is of finite representation type and not semisimple. 

Next assume that $q$ is a multiple root of $P_W(x)=0$. If $q$ is 
a multiple root of $P_{W_i}(x)=0$, for some $i$, then $\H_q(W_i)$ is 
of infinite representation type by assumption. 
Thus $\H_q(W)$ is of infinite representation type. 
If $q$ is a simple root of $P_{W_i}(x)=0$ and 
$P_{W_j}(x)=0$, for some $i\ne j$, then $\H_q(W_i)$ and $\H_q(W_j)$ are 
of finite representation type and not semisimple. By Lemma 
\ref{complexity}(3), there exist an indecomposable $\H_q(W_i)$--module 
$M_i$ and an indecomposable $\H_q(W_j)$--module $M_j$ such that 
$c_{\H_q(W_i)}(M_i)=1$ and $c_{\H_q(W_j)}(M_j)=1$. Write 
$M=M_i\otimes M_j$. 
We shall prove that the complexity of $M$ as an 
indecomposable $\H_q(W_i)\otimes\H_q(W_j)$--module is equal to $2$; 
if we use the fact that 
$c_{\H_q(W_i\times W_j)}(M)$ is the growth rate of 
${\rm Ext}^{\centerdot}(M,M)$ then the Kunneth formula implies the 
result. In a more concrete manner, 
the proof of $c_{\H_q(W_i\times W_j)}(M)=2$ is as follows. 

Let $P_i^{\centerdot}$ and $P_j^{\centerdot}$ be minimal 
projective resolutions of 
$M_i$ and $M_j$ respectively. Then 
$c_{\H_q(W_i)}(M_i)=1$ and $c_{\H_q(W_j)}(M_j)=1$ imply that there 
exists a constant $C$ such that 
$1\le {\rm dim}_K(P_i^t)\le C$ and $1\le {\rm dim}_K(P_j^t)\le C$ for all 
$t$. As $P^{\centerdot}=P_i^{\centerdot}\otimes P_j^{\centerdot}$ is a 
minimal projective resolution of $M$, we have 
\begin{equation*}
t+1\le {\rm dim}_K(P^t)=\sum_{s=0}^t {\rm dim}_K(P_i^s){\rm dim}_K(P_j^{t-s})
\le C^2(t+1). 
\end{equation*}
Therefore, the complexity of $M$ is exactly $2$. As a result, 
$\H_q(W_i)\otimes \H_q(W_j)$ is of infinite representation type 
by Lemma \ref{complexity}(2) and (3). Thus $\H_q(W)$ is also of 
infinite representation type. 
\end{proof}

\section{Type $B$ and type $D$}

To prove Theorem \ref{main theorem}, it is enough to consider 
type $B$ and type $D$ by virtue of Theorem \ref{Uno} and 
Proposition \ref{reduction}. 

Let $q$ and $Q$ be 
invertible elements of $K$.  The (two parameter) Hecke algebra
$\H_{q,Q}(B_n)$ of type  $B_n$ is the unital associative
$K$--algebra defined by generators $T_0, T_1,\dots, T_{n-1}$ and relations
$$\begin{array}{c}
(T_0+1)(T_0-Q)=0,\quad
(T_i+1)(T_i-q) =0\;(1\le i\le n-1),\\[2pt]
T_0T_1T_0T_1=T_1T_0T_1T_0,\quad
T_{i+1}T_iT_{i+1}=T_iT_{i+1}T_i\;(1\le i\le n-2),\\[2pt]
T_iT_j=T_jT_i\;(0\le i<j-1\le n-2).
\end{array}$$

The following theorem, together with Theorem \ref{Uno}, 
will allow us to assume that $-Q$ is a power of $q$.

\begin{theorem}{\cite[Theorem 4.17]{DJ}}
\label{Dipper-James}
Suppose that $Q\ne-q^f$ for any $f\in \Z$. Then 
$\H_{q,Q}(B_n)$ is Morita
equivalent to 
$\bigoplus_{m=0}^n\H_q(A_{m-1})\otimes\H_q(A_{n-m-1})$.
\end{theorem}

\begin{corollary}
\label{e=odd}
Assume that $q$ is a primitive $e^{th}$ root of unity with $e\ge 2$ as 
above. If $Q\ne-q^f$ for any $f\in\Z$ 
then $\H_{q,Q}(B_n)$ is of finite representation type 
if and only if $n<2e$.  
\end{corollary}
\begin{proof}
If $n\ge 2e$ then $\H_q(A_{n-1})$ is of infinite representation type by 
Theorem \ref{Uno}. 
Thus, $\H_{q,Q}(B_n)$ is of inifinite representation type by 
Theorem \ref{Dipper-James}. 

If $n<2e$ then one of $m$ and $n-m$ is smaller 
than $e$ for each $m$. Thus, one of $\H_q(A_{m-1})$ and $\H_q(A_{n-m-1})$ is 
semisimple and the other is of finite representation type for each $m$. 
Thus, $\H_{q,Q}(B_n)$ is of finite representation type by 
Theorem \ref{Dipper-James}. 
\end{proof}

\begin{theorem}{\cite[Theorem 1.4]{AM}}
\label{Ariki-Mathas}
Suppose that $-Q=q^f$ $(0\le f<e)$ and $e\ge 3$. 
Then $\H_{q,Q}(B_n)$ is of finite representation type if and only if
$$n<\min\{e,2\min\{f,e-f\}+4\}.$$
\end{theorem}

It is also easy to prove that Theorem \ref{Ariki-Mathas} is valid for $e=2$; 
see \cite{AM2}. 

\begin{corollary}
\label{type B}
Uno's conjecture is true if $W=W(B_n)$. 
\end{corollary}
\begin{proof}
See the argument of \cite[Introduction]{AM}. 
\end{proof}

Recall that $\H_q(D_n)$ is the $K$--algebra defined by generators 
$T_0^D,\dots,T_{n-1}^D$ and relations
$$\begin{array}{c}
(T_i^D+1)(T_i^D-q) =0\;(0\le i\le n-1),\\[2pt]
T_0^DT_2^DT_0^D=T_2^DT_0^DT_2^D,\;\; T_0^DT_i^D=T_i^DT_0^D\;(i\ne 2),\\[2pt]
T_{i+1}^DT_i^DT_{i+1}^D=T_i^DT_{i+1}^DT_i^D\;(1\le i\le n-2),\\[2pt]
T_i^DT_j^D=T_j^DT_i^D\,(1\le i<j-1\le n-2).
\end{array}$$

Now assume that $Q=1$ and denote the generators of $\H_{q,1}(B_n)$ by 
$T_0^B,\dots,T_{n-1}^B$. Then we have an algebra 
homomorphism 
\begin{equation*}
\phi: \H_q(D_n) \longrightarrow \H_{q,1}(B_n)
\end{equation*}
defined by 
$T_0^D\mapsto T_0^BT_1^BT_0^B$ and $T_i^D\mapsto T_i^B$ 
for $1\le i\le n-1$. 

$\phi$ is injective and we identify $\H_q(D_n)$ with its image. 
Then we have 
\begin{equation*}
\H_{q,1}(B_n)=\H_q(D_n)\oplus T_0^B\H_q(D_n)\text{\; and \;} 
T_0^B\H_q(D_n)=\H_q(D_n)T_0^B. 
\end{equation*}

We define an algebra automorphism $\pi$ of $\H_{q,1}(B_n)$ by 
\begin{equation*}
\pi(T_1^B)=T_0^BT_1^BT_0^B \text{\; and \;} 
\pi(T_i^B)=T_i^B\text{\; for $i\ne 1$.}
\end{equation*}
We have $\pi^2=1$ and $\pi$ induces the Dynkin automorphism of 
$\H_q(D_n)$ defined by $T_i^D\mapsto T_{1-i}^D$ for $i=0,1$ and 
$T_i^D\mapsto T_i^D$ for $i\ge 2$. 

\begin{lemma}
\label{basic lemma}
\begin{itemize}
\item[(1)] 
If $e$ is odd then $\H_{q,1}(B_n)$ is of finite representation type 
if and only if $n<2e$. 
\item[(2)] 
If $e$ is even then $\H_{q,1}(B_n)$ is of finite representation type 
if and only if $n<e$. 
\item[(3)] 
If $M$ is a semisimple $\H_{q,1}(B_n)$--module then so is 
the $\H_q(D_n)$--module ${\rm Res}(M)$. 
\end{itemize}
\end{lemma}
\begin{proof}
(1)(2) These are consequences of Corollary \ref{e=odd} 
and Theorem \ref{Ariki-Mathas}. 

\noindent
(3) We may assume that $M$ is a simple $\H_{q,1}(B_n)$--submodule without 
loss of generality. Let $N$ be a simple $\H_q(D_n)$--submodule 
of ${\rm Res}(M)$. 
Then $T_0^BN$ is also a simple $\H_q(D_n)$--submodule whose action 
is the twist of the action of $N$ by $\pi$. Since $N+T_0^BN$ is 
$\H_{q,1}(B_n)$--stable, it coincides with $M$. Hence 
${\rm Res}(M)$ is semisimple. 
\end{proof}

Let $D$ be a $\H_q(D_e)$--module which affords 
the sign representation $T_i^D\mapsto -1$, for $0\le i\le e-1$. 
We denote its projective cover by $P$. 

Recall from \cite{AM} that 
simple $\H_{q,1}(B_e)$--modules are indexed 
by Kleshchev bipartitions. Further, $\lambda=((0),(1^e))$ is 
Kleshchev and if $e\ge 3$ is even then 
the projective cover $P^{\lambda}$ of the simple 
$\H_{q,1}(B_e)$--module $D^{\lambda}$ has the property that 
$\Rad P^\lambda/\Rad^2 P^\lambda$ does not contain $D^{\lambda}$ and 
$\Rad^2 P^\lambda/\Rad^3 P^\lambda$ contains $D^{\lambda}\oplus D^{\lambda}$. 
See the proof of \cite[Theorem 4.1]{AM} in p.12. 
$D^{\lambda}$ affords the representation $T_i^B\mapsto -1$ 
$(0\le i\le n-1)$. 

\begin{lemma}
\label{pim}
Assume that $e$ is even. Then $P\simeq {\rm Res}(P^{\lambda})$. 
\end{lemma}
\begin{proof}
First note that the characteristic $l$ of the base field is odd 
since $e$ is even. 
Let $D^{\mu}$ be the simple $\H_{q,1}(B_e)$--module which affords 
the representation $T_0^B\mapsto 1$ and $T_i^B\mapsto -1$ 
$(1\le i\le e-1)$. Then ${\rm Ind}(D)=D^{\lambda}\oplus D^{\mu}$ 
since the left hand side is given by 
\begin{equation*}
T_0^B=\begin{pmatrix} 0 & 1 \\ 1 & 0 \end{pmatrix},\quad
T_i^B=\begin{pmatrix} -1 & 0 \\ 0 & -1 \end{pmatrix} \;(1\le i\le e-1),
\end{equation*}
and $T_0^B$ is diagonalizable because $l$ is odd. Now 
the surjection $P \rightarrow D$ induces 
surjective homomorphisms 
${\rm Ind}(P)\rightarrow D^{\lambda}$ and 
${\rm Ind}(P)\rightarrow D^{\mu}$. Hence these induce 
surjective homomorphisms ${\rm Ind}(P)\rightarrow P^{\lambda}$ and 
${\rm Ind}(P)\rightarrow P^{\mu}$. We have that both $P^{\lambda}$ and 
$P^{\mu}$ are direct summands of ${\rm Ind}(P)$. On the other hand, 
Mackey's formula implies that 
\begin{equation*}
{\rm Res}\left({\rm Ind}(P)\right)\simeq P\oplus\hphantom{}^{\pi}\!P, 
\end{equation*}
where $\hphantom{}^{\pi}\!P$ is the indecomposable projective 
$\H_q(D_e)$--module 
which is the twist of $P$ by $\pi$. As the twist of $D$ is $D$ itself, 
we have ${\rm Res}({\rm Ind}(P))\simeq P\oplus P$. As 
${\rm Res}(P^{\lambda})$ and ${\rm Res}(P^{\mu})$ are direct summands 
of ${\rm Res}({\rm Ind}(P))$, we conclude that 
${\rm Res}(P^{\lambda})$ and ${\rm Res}(P^{\mu})$ are isomorphic to 
$P$. 
\end{proof}

Next lemma is obvious. 

\begin{lemma}
\label{summand}
Let $A$ be a finite dimensional $K$--algebra and $B$ a $K$--subalgebra 
such that $B$ is a direct summand of $A$ as a $(B,B)$--bimodule. 
If $A$ is of finite representation type then so is $B$. 
\end{lemma}

Recall that $P_W(x)=(x^n-1)\prod_{i=1}^{n-1}\frac{x^{2i}-1}{x-1}$ in 
type $D_n$. 
Thus $q$ is a simple root of $P_W(x)=0$ if and only if 
either $e$ is odd and $e\le n< 2e$ or $e$ is even and 
$\frac{e}{2}+1\le n<e$. 

\begin{proposition}
Uno's conjecture is true if $W=W(D_n)$. 
\end{proposition}
\begin{proof}
First assume that $e$ is odd. If $n<2e$ then Lemma \ref{summand} implies 
that $\H_q(D_n)$ is of finite representation type since 
$\H_{q,1}(B_n)$ is of finite representation type 
by Lemma \ref{basic lemma}(1). If $n\ge 2e$ then it is enough to 
prove that $\H_q(D_{2e})$ is of infinite representation type by 
Lemma \ref{summand}. Using the same lemma again, we further 
know that it is enough to prove that 
$\H_q(A_{2e-1})$ is of infinite representation type. However, this 
is nothing but the result of Theorem \ref{Uno}. 

Next assume that $e$ is even. If $n<e$ then Lemma \ref{summand} implies 
that $\H_q(D_n)$ is of finite representation type since 
$\H_{q,1}(B_n)$ is of finite representation type 
by Lemma \ref{basic lemma}(2). If $n\ge e$ then it is enough to 
prove that $\H_q(D_e)$ is of infinite representation type by 
Lemma \ref{summand}. Note that $W(D_2)=W(A_1)\times W(A_1)$ and 
the conjecture is true in this case by Proposition \ref{reduction}. 
Thus we may assume that $e\ge 4$. In particular, we have $q\ne -1$, 
and this implies that ${\rm Ext}^1(D,D)=0$; if we write 
$$T_i=\begin{pmatrix} -1 & a_i \\ 0 & -1 \end{pmatrix},$$
then $T_i-q$ is invertible and thus $a_i=0$. 

Let $\overline P={\rm Res}(P^{\lambda}/\Rad^3 P^{\lambda})$. 
Lemma \ref{basic lemma}(3) and ${\rm Ext}^1(D,D)=0$ 
imply that $\overline P$ has Loewy length $3$. 
Since $\overline P$ has unique head $D$ by 
Lemma \ref{pim}, there exists a surjective 
homomorphism $P\rightarrow \overline P$. Further, as 
$\overline P$ contains $D\oplus D$ as a $\H_q(D_e)$--submodule, 
$\Rad^2 P/\Rad^3 P$ contains $D\oplus D$. On the other hand, 
${\rm Ext}^1(D,D)=0$ implies that $\Rad P/\Rad^2 P$ does not 
contain $D$. Hence we conclude that ${\rm End}_{\H_q(D_e)}(P/\Rad^3 P)$ 
is isomorphic to $K[X,Y]/(X^2,XY,Y^2)$, which is 
not isomorphic to any of the truncated polynomial rings 
$K[X]/(X^N)$ $(N=1,2,\dots)$. As we assume that the base field is a 
splitting field of $\H_q(D_e)$, this implies that $\H_q(D_e)$ is 
of infinite representation type. 
\end{proof}

\section{Appendix}

In this section, we prove Theorem \ref{q=1}. 
If the reader is familiar with the structure of the Sylow subgroups 
of exceptional Weyl groups then (s)he would not need this proof to 
know that Theorem \ref{q=1} is true. 
In the proof below, we use standard facts about the structure of 
exceptional Weyl groups; they can be found in 
\cite{Bourbaki} or \cite[2.12]{Humphreys}. 
First we consider irreducible Weyl groups. 

\bigskip
\noindent
Type $A_{n-1}$; 

$W$ has cyclic Sylow $l$--subgroups if and only if $n<2l$, and 
this is equivalent to 
the condition that $l^2$ does not divide $|W|=n!$. 

\bigskip
\noindent
Type $B_n$; 

$W$ has cyclic Sylow $l$--subgroups if and only if 
either $l>2$ and $n<2l$ or $l=2$ and $n<2$, and 
this is equivalent to 
the condition that $l^2$ does not divide $|W|=2^nn!$. 

\bigskip
\noindent
Type $D_n$; 

$W$ has cyclic Sylow $l$--subgroups if and only if either 
$l>2$ and $n<2l$ or $l=2$ and $n<2$, and this is equivalent to 
the condition $l^2$ does not divide $|W|=2^{n-1}n!$. 

\bigskip
\noindent
Type $F_4$; 

As $|W(F_4)|=2^7\cdot 3^2$, we prove that the 
Sylow $l$--subgroups for $l=2, 3$ are not cyclic. 
Let $\Delta(F_4)$ be the root system of type $F_4$. 
The long roots form a root system which is isomorphic 
to $\Delta(D_4)$. Let $\Gamma$ be the Dynkin automorphism 
group of the Dynkin diagram of type $D_4$. Then it is known 
that $W(F_4)$ is isomorphic to the semi--direct product of 
$W(D_4)$ and $\Gamma$. 

Assume that $l=2$. Since $W(D_4)$ contains $C_2\times C_2$, the 
Sylow $2$--subgroup of $W(F_4)$ cannot be cyclic. 

Assume that $l=3$. Since $\Gamma$ is isomorphic to 
the symmetric group of degree $3$, we can choose $\sigma\in\Gamma$ 
of order $3$. Let $P$ be a Sylow $3$--subgroup of $W(F_4)$ 
containing $\sigma$. As the Sylow $3$--subgroup of $W(D_4)$ is a 
cyclic group of order $3$, we have 
$|W(D_4)\cap P|\le 3$. On the other hand, as 
$\langle\sigma\rangle$ is a Sylow $3$--subgroup of 
$\Gamma$ and $|P|=9$, both $W(D_4)\cap P$ and 
$P/W(D_4)\cap P$ are isomorphic to the 
cyclic group of order $3$. Let $\tau$ be a generator of 
$W(D_4)\cap P$. Then 
we have the following split exact sequence. 
\begin{equation*}
1 \longrightarrow \langle\tau\rangle\simeq C_3 \longrightarrow 
P \longrightarrow \langle\sigma\rangle\simeq C_3 \longrightarrow 1.
\end{equation*}
As ${\rm Aut}(C_3)\simeq C_2$, $\sigma$ acts on 
$\langle\tau\rangle$ trivially and $P\simeq C_3\times C_3$. 

\bigskip
\noindent
Type $E_n$; 

Recall that $|W(E_6)|=2^7\cdot 3^4\cdot 5$, 
$|W(E_7)|=2^{10}\cdot 3^4\cdot 5\cdot 7$ and 
$|W(E_8)|=2^{14}\cdot 3^5\cdot 5^2\cdot 7$. 

Assume that $l=2$ or $l=3$. Since $W(F_4)\subset W(E_n)$, 
The Sylow $l$--subgroup of $W(E_n)$ contains $C_l\times C_l$. 
Thus it cannot be a cyclic group. 

Assume that $l=5$. Let $Q$ be the root lattice of type $E_8$ 
with scalar product normalized to $(\alpha_i,\alpha_i)=2$ for 
simple roots $\alpha_i\;(1\le i\le 8)$. Then 
$q(x)=\frac{(x,x)}{2}\;({\rm mod}\;2)$ defines a quadratic form 
on $Q/2Q\simeq {\mathbb F}_2^8$. Note that 
if we choose simple roots as a basis, we can write down 
$q(x)$ explicitly, and the computation of its Witt 
decomposition shows that 
its Witt index is $4$. Thus, by choosing a different basis, 
we may assume that $q(x)=\sum_{i=1}^4 x_{2i-1}x_{2i}$. 
Let $O_8(2)$ be the orthogonal group associated to this form. 
Then it is known that there is an exact sequence 
\begin{equation*}
1 \longrightarrow \{\pm1\} \longrightarrow W(E_8) \longrightarrow 
O_8(2) \longrightarrow 1. 
\end{equation*}

Let $q'(x)=x_1x_2+x_3^2+x_3x_4+x_4^2$ be a quadratic form on 
${\mathbb F}_2^4$. If we write $O_4^-(2)$ for the orthogonal group 
associated to this form, we know that its Sylow $5$--subgroups are 
cyclic of order $5$. 
Now explicit computation of Witt decomposition again 
shows that $q'\oplus q'$ has Witt index $4$. Thus 
$O_8(2)$ contains $C_5\times C_5$. As a result, the 
Sylow $5$--subgroup of $W(E_8)$ is isomorphic to 
$C_5\times C_5$. 

\bigskip
\noindent
Type $G_2$; 

This is the dihedral group of order $12$ and its 
Sylow $2$--subgoups are not cyclic. 

\bigskip
Now let $W$ be a general finite Weyl group. 
That is, $W$ is a product of the groups listed above. 
Then $W$ has a cyclic Sylow $l$--group 
if and only if at most one component of the product has 
a cyclic Sylow $l$--group and all the other components have 
trivial Sylow $l$--groups. This is equivalent to the condition 
that $l^2$ does not divide $|W|$. 

%\backmatter

%\include{biblio}
%\include{index}
\end{document}